\documentclass[12pt, a4paper]{article}
\usepackage{amsfonts}
\usepackage{mathrsfs}

\usepackage{latexsym}
\usepackage{graphicx}
\usepackage{amssymb}

\newtheorem{theorem}{Theorem}[section]
\newtheorem{lemma}[theorem]{Lemma}

\title{{\bf Spectral conditions for the existence of specified paths and cycles in graphs}
\thanks{Supported by the NNSF of China (No.11101057), China
Postdoctoral Science Foundation (No.20110491443) and Zhejiang
Provincial Natural Science Foundation of China(LY12A01016).}}

\author{{\bf Mingqing Zhai$^{a}$}
, {\bf Huiqiu Lin$^{b}$}\thanks{Corresponding author. E-mail addresses: mqzhai@chzu.edu.cn
(M.Zhai); huiqiulin@126.com (H.Lin); scgong@zafu.edu.cn(S.Gong).}
, {\bf Shicai Gong$^{c}$}
\\
{\footnotesize $^a$School of Mathematical Science, Chuzhou
University, Anhui, Chuzhou, 239012, China} \\
{\footnotesize $^b$Department of Mathematics, East China University of Science and Technology, Shanghai, 200237, China }\\
{\footnotesize $^c$School of Science, Zhejiang A {\&} F University,
Zhejiang, Linan, 311300, China}}

\date{}
\begin{document}
\openup 1.0\jot \maketitle

\begin{abstract}
Let $G$ be a graph with $n$ vertices and $\lambda_n(G)$ be the least
eigenvalue of its adjacency matrix of $G$. In this paper, we give
sharp bounds on the least eigenvalue of graphs without given pathes
or cycles and determine the extremal graphs. This result gives
spectral conditions for the existence of specified paths and cycles
in graphs.


\end{abstract}

\section{Introduction}

~~Let $G=(V,E)$ be the graph with vertex set $V$ and edge set $E$.
Let $|G|$ and $e(G)=|E(G)|$ be the order and the size of $G$,
respectively. A graph $G$ is said $H$-free, if $G$ does not contain
a copy of $H$. Denote by $P_t$ and $C_t$ the path and cycle of order
$t$, respectively. Let $N^d(u)=\{v|v\in V(G),d_G(v,u)=d\}$, where
$V(G)$ is the vertex set of $G$ and $d_G(v,u)$ is the distance
between $u$ and $v$. Particularly, denote $N(u)=N^1(u)$ and
$d(u)=|N(u)|$. For two disjoint vertex sets $V_1$ and $V_2$, let
$e(V_1,V_2)$ be the number of edges with one endpoint in $V_1$ and
another in $V_2$. Let $A(G)$ be the adjacency matrix of a graph $G$.
The largest modulus of an eigenvalue of $A(G)$ is called the
spectral radius of $G$ and denoted by $\rho(G)$. It is known that
for a connected graph $G$, there exists a positive unit eigenvector
corresponding to $\rho(G)$, which is called the Perron vector of
$G$.

A central problem of extremal graph theory is as follows: for a
given graph $H$, what is the maximum number of edges of an $H$-free
graph of order $n$? In the past decades, much attention has been
paid to the spectral version of above question, that is, what is the
maximum spectral radius of an $H$-free graph of order $n$? It is
known that if $H$ is a star $K_{1,r+1}$, then the extremal graph is
an $r$-regular graph (see \cite {MIN}). In recent years, Fiedler and
Nikiforov \cite{V1} solved the question for $H\cong C_n$ or $H\cong
P_n$. Nikiforov solved the case $H\cong C_{2k+1}$ for large enough
$n$ (see \cite{V2}) and the case that $H$ is a complete graph
$K_{r+1}$ or a complete bipartite graph $K_{2,r+1}$ (see \cite{V3}).
In \cite{V4}, Nikiforov showed if $G$ is a $C_4$-free graph with $n$
vertices, then $\rho^2(G)-\rho(G)-(n-1)\leq0$ with equality if and
only if $G$ is the friendship graph (now $n$ is odd). Moreover, he
conjectured if $G$ is a $C_4$-free graph with even vertices, then
$\rho^3(G)-\rho^2(G)-(n-1)\rho(G)+1\leq0$ with equality if and only
if $G$ is a star of order $n$ with $\frac n2-1$ disjoint additional
edges. Zhai and Wang \cite{Zhai} showed this conjecture is true. In
\cite {V5}, Nikiforov obtained the extremal graph for $H\cong
P_{2k+2}$ or $H\cong P_{2k+3}$ and gave a conjecture on $H\cong
C_{2k+2}$ (see \cite {V5}, Conjecture 15.). This paper solves these
two questions in the class of bipartite graphs (see Theorems
\ref{th3} and \ref{th4}).

Let $\lambda_n(G)$ be the least eigenvalue of a graph $G$ of order $n$.
It is known that $\lambda_n(G) =-\rho(G)$ for a
bipartite graph $G$ (see \cite{Doob}). Recently, researchers have
begun to pay attention to the least eigenvalues of graphs with a
given value of some well-known integer graph invariant: for
instance: order and size \cite{Bell1,Bell2,Fan1,MP}, unicyclic
graphs with a given number of pendant vertices \cite{Liu}, matching
number and independence number \cite{Fan2}, number of cut vertices
\cite{Fan3}, connectivity, chromatic number \cite{Fan4}, domination
number \cite{Zhu}. This paper also gives a spectral extremal characterization on the least eigenvalue of graphs
(see Theorems \ref{th1} and \ref{th2}).

\begin{lemma}\label{le1}
Let $G=<X,Y>$ be a bipartite graph, where $|X|\geq k$, $|Y|\geq k-1$.
If $G$ does not contain a copy of $P_{2k+1}$ with both endpoints in $X$, then
$$e(G)\leq (k-1)|X|+k|Y|-k(k-1).          \eqno(1)$$
Equality holds if and only if $G$ is isomorphic to a complete bipartite graph $K_{|X|,|Y|}$, where $|X|=k$ or $|Y|=k-1$.
\end{lemma}

\begin{lemma}\label{le2}
Let $G$ be a $P_{2k+3}$-free bipartite graph. If $G$ does not
contain a connected component isomorphic to $K_{k+1,k+1}$, then for
an arbitrary vertex $u\in V(G)$,
$$e(N(u), N^2(u))\leq (k-1)|N(u)|+k|N^2(u)|-k(k-1),       \eqno(2)$$
and equality holds if and only if the subgraph induced by $N(u)\cup N^2(u)$ is isomorphic to one of the following:\\
(i) $K_{|N(u)|,|N^2(u)|}$, where $|N(u)|=k$ or $|N^2(u)|=k-1$;\\
(ii) a graph obtained from $K_{k,|N^2(u)|}$ and an additional vertex $v\in N(u)$ by joining it with $k-1$ vertices of $N^2(u)$.
\end{lemma}

\begin{theorem}\label{th1} 
Let $G$ be a $C_t$-free graph of order $n$.\\
(i) If $t$ is odd or $n<t$, then $\lambda_n(G)\geq
-\sqrt{\lfloor\frac n2\rfloor\lceil\frac n2\rceil}$, and equality
holds if and
only if $G\cong K_{\lfloor\frac n2\rfloor, \lceil\frac n2\rceil}$.\\
(ii) If $t=2k+2$ and $n\geq t$, then $\lambda_n(G)\geq-
\sqrt{k(n-k)}$, and equality holds if and only if $G\cong
K_{k,n-k}$.
\end{theorem}

\begin{theorem}\label{th2} 
Let $G$ be a $P_t$-free graph of order $n$.\\
(i) If $n<t$, then $\lambda_n(G)\geq -\sqrt{\lfloor\frac
n2\rfloor\lceil\frac n2\rceil}$. Equality holds if and
only if $G\cong K_{\lfloor\frac n2\rfloor, \lceil\frac n2\rceil}$.\\
(ii) If $n\geq t$, then $\lambda_n(G)\geq -\sqrt{k(n-k)}$, where
$k=\lfloor\frac{t-2}2\rfloor$. Equality holds if and only if $G\cong
K_{k,n-k}$ or $G\cong K_{2,2}\cup K_1$ for $(n,t)=(5,5)$.
\end{theorem}

In the remaining part of this paper, we give the proofs of Lemmas
1.1 and 1.2, and Theorems 1.3 and 1.4.

\section{Proofs}

\noindent\textbf{ Proof of Lemma 1.1.} We use induction on $|Y|$. If
$|Y|=k-1$, then $$e(G)\leq e(K_{|X|,k-1})=(k-1)|X|.$$ So, (1)
follows and the equality holds if and only if $G\cong K_{|X|,k-1}$.

Now we suppose that the inequality holds for $|Y|\leq M$, where
$M\geq k$.

\noindent\textbf{Case 1.} There exists a vertex $u\in Y$ such that
$d_G(u)\leq k-1$, then by the hypothesis, $$e(G-u)\leq
(k-1)|X|+k(|Y|-1)-k(k-1).$$ Hence we have
$$e(G)=e(G-u)+d_G(u)<(k-1)|X|+k|Y|-k(k-1).
$$

\noindent\textbf{Case 2.} For each vertex $u\in Y$, $d_G(u)\geq k$.
We use induction on $|X|$. If $|X|=k$, then $$e(G)\leq
e(K_{k,|Y|})=k|Y|,$$ so (1) follows and the equality holds if and
only if $G\cong K_{k,|Y|}$. Suppose that the inequality holds for
$|X|\leq M$ where $M\geq k+1$.

\noindent\textbf{Subcase 2.1.} There exists a vertex $v\in X$, such
that $d_G(v)\leq k-2$, then by the hypothesis, we have $$e(G-v)\leq
(k-1)(|X|-1)+k|Y|-k(k-1).$$ Hence
$$e(G)=e(G-u)+d_G(u)<(k-1)|X|+k|Y|-k(k-1).$$

\noindent\textbf{Subcase 2.2.} For each vertex $v\in X$, $d_G(v)\geq
k-1$. We may assume that $G$ is connected. Otherwise, we use the
induction hypothesis on each component of $G$ and the inequality
follows immediately. Now we suppose that $P$ is a longest path in
$G$ with two endpoints $w_1,w_2$.

\noindent\textbf{Subcase 2.2.1.} Either $w_1$ or $w_2$ in $Y$, say
$w_1\in Y$. Then $|P\cap Y|\geq |P\cap X|$. Note that $G$ is
$P_{2k+2}$-free, then $|P\cap X|\leq k$. Since $d_G(w_1)\geq k$,
thus we have $|P\cap X|=k$ and $N_G(w_1)=P\cap X$. One can find now
each vertex in $P\cap Y$ is an endpoint of a path of the same length
with $P$. Similarly, we have $N_G(w)=P\cap X$ for any $w\in P\cap
Y$. Then $X=P\cap X$ (otherwise we can find a path with length
longer than $|P|$ which is forbidden). Thus we have $d_G(u)=k$ for
each $u\in Y$ since $d_G(u)\geq k$ and $|X|=|P\cap X|=k$. Therefore
$G\cong K_{k,|Y|}$ and then the inequality in (1) is an equality.

\noindent\textbf{Subcase 2.2.2.} Both $w_1$ and $w_2$ are in $X$.
Then $|P\cap X|=|P\cap Y|+1$. Note that $G$ does not contain a copy
of $P_{2k+1}$ with both endpoints in $X$, then $|P\cap Y|\leq k-1$.
Since $P$ is a longest path in $G$, thus $N_G(w_1)\subset P\cap Y$.
Hence we have $N_G(w_1)= P\cap Y$ and $d_G(w_1)=|P\cap Y|=k-1$ since
$d_G(w_1)\geq k-1$. Further, each vertex $u\in X$ is an endpoint of
a path of the same length with $P$. So we have $d_G(u)=k-1$ for each
vertex $u\in X$. Similarly, we have $Y=P\cap Y$. Therefore, $G\cong
K_{|X|,k-1}$ and then the inequality in (1) is an equality.

Thus we complete the proof.\ \ \ $\Box$

\noindent\textbf{Proof of Lemma 1.2.} Denote by $H$ the subgraph of
$G$ induced by $N(u)\cup N^2(u)$. By Lemma \ref{le1}, if $H$ does
not contain a copy of $P_{2k+1}$ with both endpoints in $N(u)$, then
Inequality (2) holds, and equality holds if and only if $G\cong
K_{|N(u)|,|N^2(u)|}$ with $|N(u)|=k$ or $|N^2(u)|=k-1$. Now suppose
that $H$ contains a path $P$ of order $2k+1$ with both endpoints in
$N(u)$.

Since $G$ is $P_{2k+3}$-free, $N(u)\backslash V(P)=\emptyset$ which
implies that $|N(u)|=k+1$ and $|N^2(u)|\geq k$. We distinguish the
following two cases.

\noindent\textbf{Case 1.} There is a vertex $v\in N(u)$ such that
$d_H(v)\leq k-1$. Now, $$e(H-v)\leq
e(K_{k,|N^2(u)|})=k|N^2(u)|=(k-1)|N(u)|+k|N^2(u)|-(k+1)(k-1),$$ and
equality holds if and only if $H-v\cong K_{k,|N^2(u)|}$. Therefore,
$$e(H)=e(H-v)+d_H(v)\leq (k-1)|N(u)|+k|N^2(u)|-k(k-1)$$ If equality
holds, then $d_H(v)=k-1$.

\noindent\textbf{Case 2.} For each vertex $v\in N(u)$, $d_H(v)\geq
k$. Now let $w\in N(u)$ be an endpoint of $P$. Note that $P$ is a
longest path in $H$. Then $d_H(w)\geq k$ and $N_H(w)\subseteq
V(P)\cap N^2(u)$. However, $|V(P)\cap N^2(u)|=k$. Hence
$N_H(w)=V(P)\cap N^2(u)$. Further, we can find that each vertex
$x\in N(u)$ must be an endpoint of a path of order $2k+1$ in $H$.
Correspondingly, $N_H(x)=V(P)\cap N^2(u)$. This implies that $V(P)$
induces a copy of $K_{k+1,k}$. Since $G$ is $P_{2k+3}$-free,
$N^2(u)\backslash V(P)=\emptyset$. Thus $H\cong K_{k+1,k}$. Since
$G$ does not contain a connected component isomorphic to
$K_{k+1,k+1}$, $N^3(u) \neq\emptyset$. We then have a path of order
$2k+3$ in $G$, a contradiction.

This completes the proof. \ \ \ $\Box$

%
%
%
%
%

\begin{theorem}\label{th3} 
Let $G$ be a $C_t$-free bipartite graph of order $n$.\\
(i) If $t$ is odd or $n<t$, then $\rho(G)\leq \sqrt{\lfloor\frac
n2\rfloor\lceil\frac n2\rceil}$. Equality holds if and
only if $G\cong K_{\lfloor\frac n2\rfloor, \lceil\frac n2\rceil}$.\\
(ii) If $t=2k+2$ and $n\geq t$, then $\rho(G)\leq \sqrt{k(n-k)}$.
Equality holds if and only if $G\cong K_{k,n-k}$.
\end{theorem}

\noindent\textbf{Proof.} (i) Recall that $K_{\lfloor\frac n2\rfloor,
\lceil\frac n2\rceil}$ has maximal spectral radius among all
bipartite graphs of order $n$. If $t$ is odd or $n<t$, then
$\rho(G)\leq \sqrt{\lfloor\frac n2\rfloor\lceil\frac n2\rceil}$
since $G$ is bipartite. Equality holds if and only if $G\cong
K_{\lfloor\frac n2\rfloor, \lceil\frac n2\rceil}$.

(ii) Let $G$ be a bipartite graph without $C_{2k+2}$ and
$$B=(b_{ij})_{n\times n}=A(G)^2-k(n-k)I,$$ where $I$ is the identity matrix of order $n$. Since $G$ is $C_{2k+2}$-free, $G$
does not contain a copy of $P_{2k+1}$ with both endpoints in $N(u)$
for any $u\in V(G)$. Let $b(u)$ be the sum of the $u$-th row of $B$.
Note that $(A(G)^2)_u$ is the number of walks of length 2 which are
started from $u$. Then by Lemma \ref{le1}, for any $u\in V(G)$,
\begin{eqnarray*}
b(u)&=&|N(u)|+e(N(u), N^2(u))-k(n-k)\\
&\leq& |N(u)|+(k-1)|N(u)|+k|N^2(u)|-k(k-1)-k(n-k)\\
&=& k(|N(u)|+|N^2(u)|)-k(n-1)\\
&\leq& k(n-1)-k(n-1)\\
&=& 0.
\end{eqnarray*}
If the equality holds, then by the second inequality,
$N^3(u)=\emptyset$. Further, by the first inequality and Lemma
\ref{le1},  $G-u\cong K_{|N(u)|,|N^2(u)|}$ with $|N(u)|=k$ or
$|N^2(u)|=k-1$. This implies that $G\cong K_{k,n-k}.$

Let $X$ be an
eigenvector of $G$ corresponding to $\rho(G)$ with
$\sum_{i=1}^nx_i=1$. Then
$$\rho^2(G)-k(n-k)=\sum_{i=1}^n[\rho^2(G)-k(n-k)]x_i=\sum\limits_{i=1}^n\sum\limits_{j=1}^nb_{ij}x_j
=\sum\limits_{u\in V(G)}b(u)x_u\leq 0.$$
Thus we have $\rho(G)\leq \sqrt{k(n-k)}$. If the equality holds,
then $b(u)=0$ for any $u\in V(G)$. Thus $G\cong
K_{k,n-k}$. On the contrary, if $G\cong
K_{k,n-k}$, then $\rho(G)=\sqrt{k(n-k)}.$ This completes the proof. \ \ \ $\Box$

\begin{theorem}\label{th4} 
Let $G$ be a $P_t$-free bipartite graph of order $n$.\\
(i) If $n<t$, then $\rho(G)\leq \sqrt{\lfloor\frac
n2\rfloor\lceil\frac n2\rceil}$. Equality holds if and
only if $G\cong K_{\lfloor\frac n2\rfloor, \lceil\frac n2\rceil}$.\\
(ii) If $n\geq t$, then $\rho(G)\leq \sqrt{k(n-k)}$, where
$k=\lfloor\frac{t-2}2\rfloor$. Equality holds if and only if $G\cong
K_{k,n-k}$ or $G\cong K_{2,2}\cup K_1$ for $n=t=5$.
\end{theorem}

\noindent\textbf{Proof.} Since $K_{\lfloor\frac n2\rfloor,
\lceil\frac n2\rceil}$ has maximal spectral radius among all
bipartite graphs of order $n$, (i) holds clearly. Next consider
(ii).

If $t=2k+2$, that is, $G$ is $P_{2k+2}$-free, then $G$ is also $C_{2k+2}$-free. Note that $K_{k,n-k}$ does not contain
a copy of $P_{2k+2}$. By (ii) of Theorem \ref{th3}, the result holds.

If $t=2k+3$, then $n\geq 2k+3$. Note that
$$\rho(K_{k,n-k})=\sqrt{k(n-k)}\geq \sqrt{k(k+3)}\geq (k+1)=\rho(K_{k+1,k+1}).$$
Above inequalities become equalities if and only if $n=2k+3=5$, that is, $n=t=5$. So we
may assume that $G$ does not contain a connected component
isomorphic to $K_{k+1,k+1}$. Then by Lemma \ref{le2},
$$e(N(u), N^2(u))\leq (k-1)|N(u)|+k|N^2(u)|-k(k-1)$$
for any $u\in V(G)$.

Let $X$ be an
eigenvector of $G$ corresponding to $\rho(G)$ with
$\sum_{i=1}^nx_i=1$ and $$B=(b_{ij})_{n\times n}=A^2(G)-k(n-k)I.$$ Similar to the proof of Theorem \ref{th3}, we have
$$\rho^2(G)-k(n-k)=\sum_{i=1}^n[\rho^2(G)-k(n-k)]x_i=\sum\limits_{i=1}^n\sum\limits_{j=1}^nb_{ij}x_j
=\sum\limits_{u\in V(G)}b(u)x_u\leq 0.$$
Hence, $\rho(G)\leq \sqrt{k(n-k)}$. If the equality holds,
then $b(u)=0$ for any $u\in V(G)$. This also implies that $N^3(u)=\emptyset$ and (2) becomes an equality for any $u\in V(G)$.
Note that if the case (ii) of Lemma \ref{le2} occurs, then $d_G(u)=k+1$ and $d_G(v)=k$. Thus
either $G-u\cong K_{k,|N^2(u)|}$ or $G-u\cong K_{|N(u)|,k-1}$ (otherwise, there will exist a vertex $v\in V(G)$
such that $b(v)<0$). So $G\cong
K_{k,n-k}$. we complete
the proof. \ \ \ $\Box$

The following lemma is very useful in characterizing the extremal
graph with minimal least eigenvalue.

\begin{lemma}[\cite{Cao}]\label{le3} If $G$ is a simple graph with $n$ vertices, then there
exists a spanning bipartite subgraph $H$ of $G$ such that
$\lambda_n(G)\geq \lambda_n(H)$ with equality if and only if $G\cong
H$.
\end{lemma}

Using Lemma \ref{le3} and Theorems 2.1 and 2.2, we immediately get
Theorems 1.3 and 1.4.

\noindent{\bf Acknowledgment} Thanks to Prof. Guanghui Xu who
organized a workshop in Zhejiang A {\&} F University during Aug.
04-10.

\small {

}

\begin{thebibliography}{99}
\bibitem{Bell1}F.K. Bell, D. Cvetkovi\'c, P. Rowlinson, S.K. Simi\'c, Graph
for which the least eigenvalues is minimal, I, Linear Algebra and its
Applications, 429 (2008) 234-241.

\bibitem{Bell2}F.K. Bell, D. Cvetkovi\'c, P. Rowlinson, S.K. Simi\'c, Graph
for which the least eigenvalues is minimal, II, Linear Algebra and its
Applications, 429 (2008) 2168-2179.

\bibitem{Cao}D.S. Cao, Y. Hong, The distribution of eigenvalues of graphs,
Linear Algebra and its Applications, 216 (1995) 211-224.

\bibitem{Doob}D. Cvetkovi\'{c}, M. Doob, H. Sachs, Spectra of Graphs, third ed.,
Barth, Heidelberg, 1995.

\bibitem{Fan1}Y.Z. Fan, Y. Wang, Y.B. Gao,
Minimizing the least eigenvalues of unicyclic graphs with
application to spectral spread, Linear Algebra and its
Applications, 429 (2-3) (2008) 577-588.

\bibitem{V1}M. Fiedler, V. Nikiforov, Spectral radius and Hamiltonicity of
graphs, Linear Algebra and its
Applications, 432 (2010) 2170-2173.


\bibitem{Liu}R.F. Liu, M.Q. Zhai, J.L. Shu, The least eigenvalues of unicyclic graphs
with $n$ vertices and $k$ pendant vertices, Linear Algebra and its
Applications, 431 (5-7) (2009) 657-665.

\bibitem{MIN}
H. Minc, Nonnegative matrices, John Wiley and Sons, New York, 1988.

\bibitem{V2}
V. Nikiforov, A spectral condition for odd cycles in graphs, Linear
Algebra and its Applications, 428 (2008) 1492-1498.

\bibitem{V3}
V. Nikiforov, Bounds on graph eigenvalues II, Linear Algebra and its
Applications, 427 (2007) 183-189.

\bibitem{V4}
V. Nikiforov, The maximum spectral radius of $C_4$-free graphs of
given order and size, Linear Algebra and its Applications, 430
(2009) 2898-2905.

\bibitem{V5}
V. Nikiforov, The spectral radius of graphs without paths and cycles
of specified length, Linear Algebra and its Applications, 432 (2010)
2243-2256.

\bibitem{MP}M. Petrovi\'c, B. Borovi\'canin, T. Aleksi\'c, Bicyclic
graphs for which the least eigenvalue is minimum, Linear Algebra and its
Applications, 430 (4) (2009) 1328-1335.


\bibitem{Fan2}Y.Y. Tan, Y.Z. Fan, The vertex (edge) independence number,
vertex (edge) cover number and the least eigenvalue of a graph,
Linear Algebra and its Applications,  433 (4) (2010) 790-795.

\bibitem{Fan3}Y. Wang, Y.Z. Fan, The least eigenvalue of a graph with cut
vertices, Linear Algebra and its Applications, 433 (1) (2010) 19-27.

\bibitem{Fan4}M.L. Ye, Y.Z. Fan, D. Liang, The least eigenvalue of graphs
with given connectivity, Linear Algebra and its
Applications, 430 (4) (2009)
1375-1379.


\bibitem{Zhu}B.X. Zhu, The least eigenvalue of a graph with a given domination
number, Linear Algebra and its Applications, 437 (2012) 2713-2718.


\bibitem{Zhai}M.Q. Zhai, B. Wang, Proof of a conjecture on the spectral radius of $C_4$-free
graphs, Linear Algebra and its Applications, 437 (2012) 1641-1647.




\end{thebibliography}
\end{document}